\documentstyle[twoside]{article}
\title{\bf Asymptotic expansion of Mathieu-Bessel series. II}
\author{\sc R. B. Paris\footnote{E-mail address:\ \ {\tt r.paris@abertay.ac.uk}}\\
\\
{\em Division of Computing and Mathematics,}\\
{\em Abertay University, Dundee DD1 1HG, UK}\\
}

\begin{document}
\newcommand{\bee}{\begin{equation}}
\newcommand{\ee}{\end{equation}}
\def\f#1#2{\mbox{${\textstyle \frac{#1}{#2}}$}}
\def\dfrac#1#2{\displaystyle{\frac{#1}{#2}}}
\newcommand{\fr}{\frac{1}{2}}
\newcommand{\fs}{\f{1}{2}}
\newcommand{\g}{\Gamma}
\newcommand{\om}{\omega}
\newcommand{\br}{\biggr}
\newcommand{\bl}{\biggl}
\newcommand{\ra}{\rightarrow}
\renewcommand{\topfraction}{0.9}
\renewcommand{\bottomfraction}{0.9}
\renewcommand{\textfraction}{0.05}
\newcommand{\mcol}{\multicolumn}
\newcommand{\gtwid}{\raisebox{-.8ex}{\mbox{$\stackrel{\textstyle >}{\sim}$}}}
\newcommand{\ltwid}{\raisebox{-.8ex}{\mbox{$\stackrel{\textstyle <}{\sim}$}}}
\date{}
\maketitle
\pagestyle{myheadings}
\markboth{\hfill {\it R.B. Paris} \hfill}
{\hfill {\it Asymptotics of a Mathieu-Bessel series  } \hfill}
\begin{abstract} 
We consider the asymptotic expansion of the Mathieu-Bessel series
\[S_{\nu,\gamma}^{\mu}(a,b)=\sum_{n=1}^\infty \frac{n^\gamma K_\nu(nb/a)}{(n^2+a^2)^\mu}, \qquad (\mu>0, \nu\geq 0,  b>0, \gamma\in {\bf R})\]
as $|a|\to\infty$ in $|\arg\,a|<\fs\pi$ with the other parameters held fixed, where $K_\nu(x)$ is the modified Bessel function of the second kind of order $\nu$. We employ a Mellin transform approach to determine an integral representation for $S_{\nu,\gamma}^{\mu}(a,b)$ involving the Riemann zeta function. Asymptotic evaluation of this integral involves  appropriate residue calculations. 
Numerical examples are presented to illustrate the accuracy of each type of expansion obtained. The expansion of the alternating variant of $S_{\nu,\gamma}^\mu(a,b)$ is also considered. 
\vspace{0.4cm}

\vspace{0.3cm}

\end{abstract}

\vspace{0.2cm}

\noindent $\,$\hrulefill $\,$

\vspace{0.2cm}

\begin{center}
{\bf 1. \  Introduction}
\end{center}
\setcounter{section}{1}
\setcounter{equation}{0}
\renewcommand{\theequation}{\arabic{section}.\arabic{equation}}
We consider the expansion of the Mathieu-Bessel series defined by
\bee\label{e11}
S_{\nu,\gamma}^{\mu}(a,b)=\sum_{n=1}^\infty \frac{n^\gamma K_\nu(nb/a)}{(n^2+a^2)^\mu}, \qquad (\mu>0, \nu\geq 0, b>0,\ \gamma\in {\bf R})\ee
for large complex values of the parameter $a$,
where $K_\nu(x)$ is the modified Bessel function of the second kind. 
This is a sequel to the paper \cite{P19} in which the series with the modified Bessel function replaced by the Bessel functions $J_\nu(x)$ and $Y_\nu(x)$ were discussed for $a\to+\infty$. The case $\mu=0$ and large values of $a$ is discussed in \cite{P18}.

The functional series (without the Bessel function) with $\gamma=1$, $\mu=2$ was first introduced by Mathieu in his 1890 book \cite{EM} dealing with the elasticity of solid bodies. Much effort has been devoted to the determination of bounds and integral representations for the `algebraic' Mathieu series (with $\mu>1$) when the parameter $a>0$; see \cite{PST} and the references therein. 
Recent work of 
Gerhold and Tomovski \cite{GT} has extended the asymptotic study of algebraic Mathieu series by considering the large-$a$ expansions of the trigonometric Mathieu series (with $x$ held fixed)
\bee\label{e100}
\sum_{n=1}^\infty\frac{n \sin nx}{(n^2+a^2)^\mu},\quad\mbox{and}\quad \sum_{n=1}^\infty\frac{n \cos nx}{(n^2+a^2)^\mu}.
\ee
The coefficients in these expansions involve the polylogarithm function $\mbox{Li}_\alpha(z)=\sum_{n\geq 1} n^{-\alpha} z^n$, with $z=e^{ix}$. Earlier asymptotic results on the algebraic Mathieu series have been given in \cite{ P13, VZ}; a special case when $\gamma$ is an even integer results in a finite algebraic expansion together with an exponentially small contribution \cite{P16}. 

We apply a Mellin transform approach to determine the large-$a$ asymptotics of the sum in (\ref{e11}) and its alternating variant. This reduces the problem to a residue evaluation associated with an integral representation of $S_{\nu,\gamma}^\mu(a,b)$ in the complex plane that involves the Riemann zeta function $\zeta(s)$. The exact nature of the pole structure of the integrand is found to depend sensitively on the various parameters.
In the application of the Mellin transform method we shall require the following estimates for the gamma function and the Riemann zeta function. For real $\sigma$ and $t$, we have the estimates
\bee\label{e13}
\g(\sigma\pm it)=O(t^{\sigma-\frac{1}{2}}e^{-\frac{1}{2}\pi t}),\qquad |\zeta(\sigma\pm it)|=O(t^{\Omega(\sigma)} \log^\alpha t)\quad (t\ra+\infty),
\ee
where $\Omega(\sigma)=0$ ($\sigma>1$), $\fs-\fs\sigma$ ($0\leq\sigma\leq 1$), $\fs-\sigma$ ($\sigma<0$) and $\alpha=1$ ($0\leq\sigma\leq 1$), $\alpha=0$ otherwise \cite[p.~95]{ECT}. The zeta function $\zeta(s)$ has a simple pole of unit residue at $s=1$, $\zeta(0)=-\fs$ and has trivial zeros at $s=-2k$, $k=1, 2, \ldots\ $.
In the neighbourhood of $s=1$ we have the expansion
\bee\label{e15}
\zeta(1+\epsilon)=\frac{1}{\epsilon}\{1+\gamma_0\epsilon-\gamma_1\epsilon^2+O(\epsilon^3)\}\qquad (\epsilon\to0),
\ee
where $\gamma_0=0.5772156\ldots$ is the Euler-Mascheroni constant and $\gamma_1=-0.0728158\ldots$ is the first Stieltjes constant. The expansion of the Pochhammer symbol $(\alpha)_r=\g(\alpha+r)/\g(\alpha)$ in the neighbourhood of $\alpha$ is
\bee\label{e16}
(\alpha+\epsilon)_r=(\alpha)_r\bl\{1+\epsilon\Delta\psi(\alpha+r)+\fs\epsilon^2(\Delta^2\psi(\alpha+r)+\Delta\psi'(\alpha+r))+O(\epsilon^3)\br\}
\ee
as $\epsilon\to0$, where
\[\Delta\psi(\alpha+r):=\psi(\alpha+r)-\psi(\alpha)\]
and $\psi$ is the logarithmic derivative of the gamma function. An analogous definition applies when $\psi$ is replaced by the derivative $\psi'$.

\vspace{0.6cm}

\begin{center}
{\bf 2. \ The asymptotic expansion of $S_{\nu,\gamma}^\mu(a,b)$ for $|a|\to\infty$}
\end{center}
\setcounter{section}{2}
\setcounter{equation}{0}
\renewcommand{\theequation}{\arabic{section}.\arabic{equation}}
Let $\mu>0$, $\nu\geq0$, $b>0$ and $\gamma$ be a real parameter with $a$ complex-valued satisfying $|\arg\,a|<\fs\pi$.
We consider the asymptotic expansion of the Mathieu-Bessel series 
\bee\label{e21}
S_{\nu,\gamma}^{\mu}(a,b)=\sum_{n=1}^\infty \frac{n^\gamma K_\nu(nb/a)}{(n^2+a^2)^\mu}
\ee
as $|a|\to\infty$ in $|\arg\,a|<\fs\pi$, where $K_\nu(x)$ is the modified Bessel function. Following the procedure described in \cite{P19},
we have the Mellin transform of the function $x^\gamma K_\nu(bx)/(1+x^2)^\mu$ given by \cite{M}
\bee\label{e22a}
H(s)\equiv H_{\nu,\gamma}^\mu(s):=\int_0^\infty \frac{x^{\gamma+s-1}K_\nu(bx)}{(1+x^2)^\mu}\,dx=H_1(s,\nu)+H_2(s,\nu)+H_2(s,-\nu),
\ee
where\footnote{For ease of presentation we do not display the parameter dependence of $H_1(s,\nu)$ and $H_2(s,\nu)$ on $\mu$ and $\gamma$.}
\begin{eqnarray*}
H_1(s,\nu)&=&\frac{\pi^2 \chi^{\mu-\gamma/2-s/2}}{4\sin \pi(\lambda_\nu\!-\!\mu) \sin \pi(\lambda_{-\nu}\!-\!\mu)}{}_1{\bf F}_2\bl(\!\!\begin{array}{c}\mu\\1\!+\!\mu\!-\!\lambda_\nu,\,1\!+\!\mu\!-\!\lambda_{-\nu}\end{array}\bl|\,-\chi\br),\\
H_2(s,\nu)&=&\frac{\pi^2 \chi^{\nu/2} \g(\lambda_\nu)}{4\g(\mu) \sin \pi\nu \sin \pi(\lambda_\nu\!-\!\mu)}
{}_1{\bf F}_2\bl(\!\!\begin{array}{c}\ \ \lambda_\nu\\1\!-\!\mu\!+\!\lambda_\nu,\,1\!+\!\nu\end{array}\bl|\,-\chi\br),\\
%H_3(s)&=&-\frac{\pi^2 \chi^{-\nu/2} \g(\lambda_{\!-})}{4\g(\mu) \sin \pi\nu \sin \pi(\lambda_{\!-}\!-\!\mu)}
%{}_1{\bf F}_2\bl(\!\!\begin{array}{c}\\ \lambda_{-}\\1\!-\!\mu\!+\!\lambda_{-},1\!-\!\nu\end{array}\bl|\,-\chi\br)
\end{eqnarray*}
valid in $\Re (s)>-\gamma+\nu$,
with $\nu$ non-integer and
\[\lambda_\nu\equiv\lambda_\nu(s):=\fs(\gamma+s+\nu),\qquad \chi:=\f{1}{4}b^2.\]
The hypergeometric functions appearing in $H(s)$ are defined by 
\[{}_1{\bf F}_2\bl(\!\!\begin{array}{c}\alpha\\ \beta, \gamma\end{array}\bl|z\br)=\frac{1}{\g(\beta) \g(\gamma)}\,{}_1 F_2\bl(\!\!\begin{array}{c}\alpha\\ \beta, \gamma\end{array}\bl|z\br)=\frac{1}{\g(\beta) \g(\gamma)}\,\sum_{n=0}^\infty \frac{(\alpha)_n}{(\beta)_n (\gamma)_n}\,\frac{z^n}{n!}.\]
%where $(\alpha)_n=\g(\alpha+n)/\g(\alpha)$ is the Pochhammer symbol.
Then we obtain the integral representation upon use of the Mellin inversion theorem (see, for example, \cite[p.~118]{PK})
\bee\label{e23}
S_{\nu,\gamma}^\mu(a,b)=
\frac{a^{\gamma-2\mu}}{2\pi i} \int_{c-\infty i}^{c+\infty i}H(s)\zeta(s)a^s\,ds,
\ee
where $c>\max\{-\gamma+\nu,1\}$ and $\zeta(s)$ is the Riemann zeta function

We now need to examine the pole structure of $H(s)$. The function ${}_1{\bf F}_2(z)$ is an entire function of its parameters (and hence of the variable $s$) \cite[(16.2.5)]{DLMF} and consequently has no poles. The poles of $H(s)$ are situated at the two sequences (when $\nu$ is non-integer) $\lambda_{\pm\nu}=-k$ and the two sequences $\lambda_{\pm\nu}-\mu=\pm k$, $k=0, 1, 2, \ldots\ $. In Appendix A, it is shown that $H(s)$ is regular at the points $\lambda_{\pm\nu}-\mu=\pm k$ due to a cancellation of terms. Thus, the only singularities of $H(s)$ are those corresponding to $\lambda_{\pm\nu}=-k$; that is, at the two infinite sequences $s=-\gamma+\nu-2k$ and $s=-\gamma-\nu-2k$, $k=0, 1, 2 \ldots\ $. In addition, there is a pole of the integrand in (\ref{e23}) at $s=1$ resulting from $\zeta(s)$.

We consider the integral in (\ref{e23}) taken round the rectangular contour with vertices at $c\pm iT$ and $-c'\pm iT$,
where $c'>0$ and $T>0$. Then the leading behaviour for large $|s|$ of the hypergeometric function in $H_1(s,\nu)$ is $O(1)/(\g(1-\lambda_\nu+\mu)\g(1-\lambda_{-\nu}+\mu))$, whereas that of the hypergeometric function in $H_2(s,\nu)$ is $\chi^{-\nu/2}J_{\nu}(b)/\g(1+\lambda_\nu-\mu)$, where $J_\nu(b)$ is the Bessel function of the first kind.
We therefore find on the upper and lower sides of the rectangle $s=\sigma\pm iT$, $-c'\leq\sigma\leq c$, with $\lambda_\nu=\fs(\sigma+\gamma+\nu)\pm\fs iT$, that
\[|H(s)|=O(T^{\gamma+\sigma-2\mu-1}e^{-\pi T/2})+O(T^{\mu-1}e^{-\pi T/2})=O(T^{\alpha-1}e^{-\pi T/2})\]
as $T\to\infty$ by the first formula in (\ref{e13}), where $\alpha:=\max\{\gamma+\sigma-2\mu, \mu\}$.
Thus, the modulus of the integrand on these paths has the order estimate 
\[|H(s)a^s \zeta(s)|=O(T^{\Omega(\sigma)+\alpha-1} \log\,T e^{-\Delta T}),\qquad\Delta=\fs\pi-|\arg\,a|,\]
where $\Omega(\sigma)$ is specified in (\ref{e13}). Provided $|\arg\,a|<\fs\pi$, the contribution from the upper and lower sides of the rectangle therefore vanishes as $T\to\infty$.

Displacement of the integration path to the left over the simple poles at $s=1$ and the two sequences at $s=-(\gamma\pm\nu+2k)$ (when $\gamma\pm\nu$ is not a negative odd integer and $\nu$ is non-integer) then yields the following result.
\newtheorem{theorem}{Theorem}
\begin{theorem}$\!\!\!.$ \ Let  $\chi=b^2/4$ and $\gamma\pm\nu\neq -1, -3, \ldots\ $ with $\nu>0$ non-integer. Then we have the asymptotic expansion
\bee\label{e24}
S_{\nu,\gamma}^\mu(a,b)-a^{\gamma-2\mu+1} H_{\nu,\gamma}^\mu(1)\sim R(a;\nu)+R(a;-\nu)
\ee
as $|a|\to\infty$ in $|\arg\,a|<\fs\pi$, where 
\bee\label{e24c}
R(a;\nu):=\frac{1}{2}a^{-\nu-2\mu}\chi^{\nu/2}\g(-\nu) \sum_{k=0}^\infty\frac{(-)^k}{k!}\,(\mu)_k\zeta(-\gamma\!-\!\nu\!-\!2k) F_k^{(\mu)}(\nu;\chi)\,a^{-2k}
\ee
and
\bee\label{e24a}
F_k^{(\mu)}(\nu;\chi):={}_1F_2\bl(\!\!\begin{array}{c}-k\\ 1\!-\!\mu\!-\!k, 1\!+\!\nu\end{array}\bl|\,-\chi\br)=
\sum_{r=0}^k\frac{(-k)_r (-\chi)^r}{(1+\nu)_r (1-\mu-k)_r r!}~.
\ee 
\end{theorem}

The $F_k^{(\mu)}(\nu;\chi)$ are polynomials in $\chi$ of degree $k$.
Note that in the special case $\mu=1$ these polynomials are given by
\bee\label{e25}
F^{(1)}_k(\nu;\chi)=\sum_{r=0}^k\frac{(-\chi)^r}{(1+\nu)_r r!}.
\ee
As $k\to\infty$ the limiting value of $F^{(1)}_k(\nu;\chi)$ is $\g(1+\nu) \chi^{-\nu/2}J_\nu(b)$.

The situation covered by Theorem 1 corresponds to the case when the poles of the integrand in (\ref{e23}) at the two sequences  $s=-\gamma\pm\nu-2k$ and $s=1$ are all simple; that is, when $\nu$ is non-integer and $\gamma\pm\nu\neq -1, -3, \ldots\ $. For certain values of the parameters $\mu$, $\nu$ and $\gamma$, however, higher-order poles can arise and the residues need to be evaluated accordingly. For example, if $\gamma\pm\nu$ is a negative odd integer the pole at $s=1$ becomes a double pole. When $\nu$ assumes integer values all but a finite number of the poles of the sequence $s=-\gamma+\nu-2k$ are double. If, in addition, $\gamma\pm\nu$ is a negative odd integer, there are simple and double poles with a treble pole at $s=1$. 

We do not examine all these cases in generality, but confine attention in the next sections to some illustrative examples. In several places it will be necessary to determine values of the ${}_1F_2$ hypergeometric function when some of its parameters
are perturbed by a small amount $\epsilon$. We shall denote the $O(\epsilon)$, $O(\epsilon^2)$ contributions to the unperturbed function by the generic quantities $F_*^{(1)}$ and $F_*^{(2)}$, respectively. It is important to stress that the precise definition of these quantities will, of course, vary according to the particular situation under consideration.
\vspace{0.6cm}

\begin{center}
{\bf 3. \ Special cases}
\end{center}
\setcounter{section}{3}
\setcounter{equation}{0}
\renewcommand{\theequation}{\arabic{section}.\arabic{equation}}
In this section we shall examine the expansions of $S_{\nu,\gamma}^\mu(a,b)$ when $\nu=0, 1$ and when $\gamma+\nu$ is a negative odd integer with $\nu$ non-integer. In the first case we need to determine the limiting value of the asymptotic sums on the right-hand side of (\ref{e24}) and in the second case we have to deal with a double pole at $s=1$.
\vspace{0.4cm}

\noindent{\bf 3.1\ Integer values of $\nu$.}\ \ We present here only the details of the cases $\nu=0$ and $1$; other integer values of $\nu$ can be treated similarly. Alternatively, it is possible to use the recurrence relation
\[K_{n+2}(x)=K_n(x)+\frac{2(n+1)}{x} K_{n+1}(x)\qquad (n=0, 1, 2, \ldots)\]
to generate the expansion for $\nu=2, 3, \ldots\,$. 

When $\nu=0$ ($\gamma\neq -1, -3, \ldots$) there is a simple pole at $s=1$ and an infinite sequence of double poles at $s=-\gamma-2k$, $k=0, 1, 2 \ldots\ $. Setting $\nu=\epsilon$, $\epsilon\to0$, we have
\[R(a;\pm\epsilon)=\mp\frac{a^{-2\mu}}{2\epsilon}\bl\{1\mp\epsilon\bl(\log\,\frac{2a}{b}-\gamma_0\br)+\cdots\br\}\hspace{6cm}\]
\[\hspace{4cm}\times\sum_{k=0}^\infty\frac{(-)^k}{k!} (\mu)_k \zeta(-\gamma-2k) F_k^{(\mu)}(0;\chi)a^{-2k}\bl\{1\mp\epsilon\frac{\zeta'(-\gamma-2k)}{\zeta(-\gamma-2k)}\mp \epsilon F_*^{(1)}+\cdots\br\},\]
where $\gamma_0$ is the Euler-Mascheroni constant and $F_*^{(1)}$ is the $O(\epsilon)$ contribution to $F_k^{(\mu)}(\epsilon;\chi)$. This latter quantity can be determined by means of the following lemma.
\newtheorem{lemma}{Lemma}
\begin{lemma}$\!\!\!.$ \ Let $\alpha$, $\beta$, $\gamma$ be arbitrary complex parameters and $\epsilon_1$, $\epsilon_2$, $\epsilon_3$ be small quantities, not all of which are identically zero. If $\epsilon_1\neq0$, it is assumed that $\alpha\neq 0, -1, -2, \ldots\,$. Then as $\epsilon_1$, $\epsilon_2$, $\epsilon_3\to 0$ we have
\[{}_1F_2\bl(\!\!\begin{array}{c}\alpha+\epsilon_1\\\beta+\epsilon_2, \gamma+\epsilon_3\end{array}\bl|z\br)=
{}_1F_2\bl(\!\!\begin{array}{c}\alpha\\\beta, \gamma\end{array}\bl|z\br)+\sum_{r=1}^\infty\frac{(\alpha)_r z^r}{(\beta)_r(\gamma)_r r!}d_r(\epsilon_1,\epsilon_2,\epsilon_3)+O(\epsilon_1^2,\epsilon_2^2,\epsilon_3^2),\]
where
\[d_r(\epsilon_1,\epsilon_2,\epsilon_3):=\epsilon_1\Delta\psi(\alpha+r)-\epsilon_2\Delta\psi(\beta+r)-\epsilon_3\Delta\psi(\gamma+r)\]
and
\[\Delta\psi(\alpha+r):=\psi(\alpha+r)-\psi(\alpha),\]
with $\psi(x)$ being the logarithmic derivative of the gamma function.
\end{lemma}
\noindent {\it Proof.}\ \ 
Use of the result (\ref{e16})
in the series expansion of the ${}_1F_2(z)$ function on the left-hand side for each of the three parameters, together with the fact that $\Delta\psi(\alpha)=0$, immediately yields the stated result.\hfill$\Box$
\vspace{0.2cm}

It then follows that (with $\epsilon_1=\epsilon_2=0$ in Lemma 1)
\[F_k^{(\mu)}(\epsilon;\chi)={}_1F_2\bl(\!\!\begin{array}{c}-k\\1-\mu-k,1+\epsilon\end{array}\bl|\,-\chi\br)=F_k^{(\mu)}(0;\chi)\{1-\epsilon F_*^{(1)}+O(\epsilon^2)\},\]
where
\bee\label{e31}
F_*^{(1)}=\frac{1}{F_k^{(\mu)}(0;\chi)}\sum_{r=1}^k\frac{(-k)_r (-\chi)^r}{(1-k-\mu)_r (r!)^2}\,\Delta\psi(1+r).
\ee
The $O(\epsilon^{-1})$ terms present in the sum $R(a;\epsilon)+R(a;-\epsilon)$ on the right-hand side of (\ref{e24}) 
cancel, thereby leading to the residue contribution given by
\bee\label{e31a}
a^{-2\mu}\sum_{k=0}^\infty \frac{(-)^k}{k!} (\mu)_k \zeta(-\gamma-2k) F_k^{(\mu)}(0;\chi) a^{-2k} \bl\{\log\,\frac{2a}{b}-\gamma_0+\frac{\zeta'(-\gamma-2k)}{\zeta(-\gamma-2k)}+F_*^{(1)}\br\}.
\ee

Then we have the following theorem:
\begin{theorem}$\!\!\!.$ \ When $\nu=0$ and $\gamma\neq -1, -3, \ldots$ we have the expansion
\[S^\mu_{0,\gamma}(a,b)-a^{\gamma-2\mu+1}H_{0,\gamma}^\mu(1) \hspace{10cm}\]
\bee\label{e32}
\sim a^{-2\mu}\sum_{k=0}^\infty \frac{(-)^k}{k!} (\mu)_k \zeta(-\gamma-2k) F_k^{(\mu)}(0;\chi) a^{-2k} \bl\{\log\,\frac{2a}{b}-\gamma_0+\frac{\zeta'(-\gamma-2k)}{\zeta(-\gamma-2k)}+F_*^{(1)}\br\}
\ee
as $|a|\to\infty$ in $|\arg\,a|<\fs\pi$, where $F_k^{(\mu)}(0;\chi)$ and $F_*^{(1)}$ are defined in (\ref{e24a}) and (\ref{e31}).
\end{theorem}

If $\gamma=2m$, $m=0, 1, 2, \ldots\,$, the above expansion simplifies on account of the trivial zeros of the zeta function and the fact that $\zeta(0)=-\fs$ to yield
\[S^\mu_{0,2m}(a,b)-a^{2m-2\mu+1}H_{0,2m}^\mu(1)+\frac{1}{2}\bl(\log\,\frac{2a}{b}-\gamma_0\br)\delta_{m,0}\hspace{5cm}\]
\bee\label{e33}
\hspace{5cm}\sim a^{-2\mu}\sum_{k=0}^\infty \frac{(-)^k}{k!} (\mu)_k \zeta'(-2m-2k)F_k^{(\mu)}(0;\chi)a^{-2k}
\ee
as $|a|\to\infty$ in $|\arg\,a|<\fs\pi$,
where $\delta_{m,0}$ is the Kronecker symbol.

In the case $\nu=1$, we have simple poles at $s=1$ and $s=1-\gamma$ ($\gamma\neq0$) and double poles at $s=-\gamma-1-2k$, $k=0, 1, 2, \ldots\,$. A similar limiting procedure then yields the expansion (see Appendix B for details):
\begin{theorem}$\!\!\!.$ \ When $\nu=1$ and $\gamma\neq 0, -2, -4, \ldots$ we have the expansion
\[S^\mu_{1,\gamma}(a,b)-a^{\gamma-2\mu+1}H_{1,\gamma}^\mu(1)\sim\frac{1}{2}a^{1-2\mu} \chi^{-1/2} \sum_{k=0}^\infty \frac{(-)^k}{k!} (\mu)_k \zeta(1-\gamma-2k) a^{-2k}\]
\bee\label{e34}
+a^{-1-2\mu}\chi^{1/2}\sum_{k=0}^\infty\frac{(-)^{k+1}}{k!} (\mu)_k \zeta(-\gamma\!-\!1\!-\!2k) F_k^{(\mu)}(1;\chi) a^{-2k}
\bl\{\log\,\frac{2a}{b}-\gamma_0+\frac{1}{2}+\frac{\zeta'(-\gamma\!-\!1\!-\!2k)}{\zeta(-\gamma\!-\!1\!-\!2k)}+F_*^{(1)}\br\}
\ee
as $|a|\to\infty$ in $|\arg\,a|<\fs\pi$, where $F_k^{(\mu)}(1;\chi)$ is defined in (\ref{e24a}) and
\[F_*^{(1)}=\frac{1}{F_k^{(\mu)}(1;\chi)}\sum_{r=1}^k\frac{(-k)_r (-\chi)^r}{(2)_r(1-k-\mu)_r r!} \bl(\Delta\psi(1+r)-\frac{r}{2(1+r)}\br)\]
with $\chi=b^2/4$.
\end{theorem}

If $\gamma=2m+1$, $m=0, 1, 2, \ldots\,$, the expansion (\ref{e34}) simplifies to yield
\[S^\mu_{1,2m+1}(a,b)-a^{2m-2\mu+2}H_{1,2m+1}^\mu(1)+\frac{1}{4}a^{1-2\mu}\chi^{-1/2} \delta_{m,0}\hspace{6cm}\]
\bee\label{e35}
\hspace{4cm}\sim a^{-1-2\mu} \chi^{1/2}\sum_{k=0}^\infty\frac{(-)^{k+1}}{k!} (\mu)_k \zeta'(-2m\!-\!2\!-\!2k) F_k^{(\mu)}(1;\chi) a^{-2k}
\ee
as $|a|\to\infty$ in $|\arg\,a|<\fs\pi$.
\vspace{0.2cm}

\noindent{\bf Remark 3.1}\ \ The value of $H_{\nu,\gamma}^\mu(1)$ for $\nu=0, 1$ can be obtained by numerical integration of the integral in (\ref{e22a}) when $\gamma>\nu-1$; otherwise a limiting procedure has to be employed to deal with $H_2(s,\nu)+H_2(s,-\nu)$. Alternatively, we can use the fact that \cite{M}
\[\int_0^\infty \frac{x^\gamma K_n(bx)}{(1+x^2)^\mu}\,dx=\frac{1}{4\g(\mu)}\,G^{3,1}_{1,3}\bl(\begin{array}{c}\fs\!-\!\fs\gamma\\0, 0, \mu\!-\!\fs\!-\!\fs\gamma\end{array}\bl|\chi\br)\qquad (n=0, 1, 2, \ldots).\]
The Meijer $G$-function on the right-hand side supplies the analytic continuation of the integral into $\gamma\leq n-1$.
\vspace{0.4cm}

\noindent{\bf 3.2\ The case $\gamma+\nu=-(2m+1)$, $\nu$ non-integer.}\ \ When $\gamma+\nu=-(2m+1)$, $m=0, 1, 2, \ldots\,$ and $\nu$ is non-integer, there is a double pole at $s=1$ with all other poles being simple. Then, from (\ref{e24}), we have
\[S_{\nu,\gamma}^\mu(a,b)- \mbox{Res}\,(1)\sim\mathop{R(a;\nu)}_{\scriptstyle k\neq m}+R(a;-\nu),\]
where $\mbox{Res}\,(1)$ is the residue at $s=1$ and in the asymptotic sum $R(a;\nu)$ in (\ref{e24c}) the term corresponding to $k=m$ is omitted.

To determine the residue at $s=1$, we first write the function $H(s)$ in (\ref{e22a}) as
\bee\label{e36}
H(s)={\hat H}(s)+H_2(s,\nu),\qquad {\hat H}(s):=H_1(s,\nu)+H_2(s,-\nu),%\frac{\chi^{\nu/2} \g(-\nu)}{4\g(\mu)} \g(\fs 
\ee
where ${\hat H}(s)$ is regular at $s=1$.
The residue of $H(s) \zeta(s) a^{s}$ at $s=1$ now consists of two parts: one with a simple pole with residue $a{\hat H}(1)$, and the other with a double pole resulting from the factor $\g(\fs s-m-\fs)$ present in $H_2(s,\nu)$ and $\zeta(s)$.

The double-pole contribution can be determined by putting $s=1+\epsilon$, $\epsilon\to0$ to find
\[H_2(1+\epsilon,\nu)=\frac{\chi^{\nu/2} \g(-\nu)}{4\g(\mu)}a^{1+\epsilon}\g(-m+\fs\epsilon)\g(\mu+m-\fs\epsilon)\zeta(1+\epsilon)\hspace{3cm}\]
\bee\label{e36a}
\hspace{5cm}\times{}_1F_2\bl(\begin{array}{c}\!\!\!-m+\fs \epsilon\\\,1-\mu-m+\fs\epsilon, 1+\nu\end{array}\bl| -\chi\br).
\ee
The treatment of the $O(\epsilon)$ contribution of the above ${}_1F_2$ function is not covered by Lemma 1 (since $\alpha=-m$) and is considered in Appendix C. We find from (\ref{c3}) and (\ref{c4}) (with $n=m+1$ and $\alpha=\mu-1$) together with use of the expansion in (\ref{e15}) that
\[H_2(1+\epsilon,\nu)=\frac{(-)^m a\chi^{\nu/2} \g(-\nu)}{2\g(\mu)m!}\,\frac{F_m^{(\mu)}(\nu;\chi)}{\epsilon^2}\bl\{1+\epsilon(\log\,a+\gamma_0+\frac{1}{2}\psi(m+1)-\frac{1}{2}\psi(m+\mu)+\frac{1}{2}F^{(1)}_*)+ O(\epsilon^2)\br\},\]
where $F_m^{(\mu)}(\nu;\chi)$ is given by (\ref{e24a}) and, provided $\mu\neq 1, 2, \ldots\,$,
\[F^{(1)}_*=\frac{1}{F_m^{(\mu)}(\nu;\chi)}\bl\{\sum_{r=0}^m\frac{(-m)_r (-\chi)^r}{(1+\nu)_r(1-m-\mu)_r r!}\{\Delta\psi(m+1-r)-\Delta\psi(m+\mu-r\}\hspace{2cm}\]
\[\hspace{3cm}-\frac{(-\chi)^{m+1}}{(m+1)(1+\nu)_{m+1}(\mu-1)_{m+1}}\,{}_2F_3\bl(\begin{array}{c} 1,1\\m+2,m+2+\nu,2-\mu\end{array}\bl| -\chi\br)\br\}.\]
Thus, when $\gamma+\nu=-(2m+1)$, $\nu$ non-integer, we obtain
\[\mbox{Res}\,(1)=a^{\gamma-2\mu+1}\bl\{ {\hat H}(1)+
\frac{(-)^m\chi^{\nu/2} \g(-\nu)}{4\g(\mu)m!}\,F_m^{(\mu)}(\nu;\chi)\hspace{5cm}\]
\bee\label{e37}
\hspace{4cm}\times\{2\log\,a+2\gamma_0+\psi(m+1)-\psi(m+\mu)+F^{(1)}_*\}\br\}.
\ee
Then we obtain the result:
\begin{theorem}$\!\!\!.$ \ When $\gamma+\nu=-(2m+1)$, $m=0, 1, 2, \ldots$, $\mu\neq 1, 2, \ldots$ and $\nu\neq 0, 1, 2, \ldots$ we have the expansion
\bee\label{e38}
S_{\nu,\gamma}^\mu(a,b)- \mbox{Res}\,(1)\sim\mathop{R(a;\nu)}_{\scriptstyle k\neq m}+R(a;-\nu)
\ee
as $|a|\to\infty$ in $|\arg\,a|<\fs\pi$. The quantity
$\mbox{Res}\,(1)$ is given in (\ref{e37}) and in the asymptotic sum $R(a;\nu)$ in (\ref{e24c}) the term corresponding to $k=m$ is omitted.
A similar result with $\nu\to -\nu$ applies when $\gamma-\nu=-(2m+1)$.
\end{theorem}
\vspace{0.6cm}

\begin{center}
{\bf 4. \ The case $\gamma+\nu=-1$, $\mu=1$}
\end{center}
\setcounter{section}{4}
\setcounter{equation}{0}
\renewcommand{\theequation}{\arabic{section}.\arabic{equation}}
The case of Theorem 4 when $\mu=1, 2, \ldots$ and $\gamma+\nu=-2m-1$ requires a separate discussion since the poles at $s=2m+1-2k$ and the {\it apparent\/} singularities at $s=2m+1+2\mu\pm2k$, $k=0,1,2, \ldots\,$, coincide for an infinite number of $k$ values.
We present the details of the specific case $\gamma+\nu=-1$, $\mu=1$, with first $\nu$ non-integer and then a case when $\nu$ is an integer, and consider the sum
\[S^1_{\nu,\gamma}(a,b)=\sum_{n=1}^\infty \frac{n^{-\nu-1} K_\nu(nb/a)}{n^2+a^2}.\]
\vspace{0.4cm}

\noindent{\bf 4.1\ Non-integer $\nu$.}\ \ 
The functions $H_2(s,\pm\nu)$ in (\ref{e22a}) become
\[H_2(s,\pm\nu)=\mp\frac{\pi^2}{4\sin \pi\nu}\,\frac{J_{\pm\nu}(b)}{\sin \fs\pi(s+\gamma\pm\nu)}.\]
The poles of $H_2(s,-\nu)$ are at $s=1+2\nu-2k$, $k=0, 1, 2, \ldots$ and produce the residue contribution to the integral (\ref{e23}) given by $R(a;-\nu)$ defined in (\ref{e24c}). The component of $H(s)$ given by ${\hat H}(s):=H_1(s,\nu)+H_2(s,\nu)$ is
\bee\label{e41}
{\hat H}(s)=\frac{1}{4}\chi^{(3-s+\nu)/2} \g\bl(\frac{s-3}{2}\br) \g\bl(\frac{s-3}{2}-\nu\br)\,{}_1F_2\bl(\!\!\begin{array}{c}1\\\f{5}{2}\!-\!\fs s, \f{5}{2}\!-\!\fs s\!+\!\nu\end{array}\bl|-\chi\br)-\frac{\pi^2}{4\sin \pi\nu}\,\frac{J_{\nu}(b)}{\sin \fs\pi(s\!-\!1)}.
\ee

The residue of ${\hat H}(s)$ at the poles $s=1-2k$, $k=1, 2, \ldots$ is easily seen to be
\[\frac{(-)^{k+1}}{2} \chi^{k+1+\nu/2}\,\frac{\g(-k-1-\nu)}{(k+1)!}\,{}_1F_2\bl(\!\!\begin{array}{c}1\\k+2, k+2+\nu\end{array}\bl|-\chi\br)-\frac{(-)^k\pi J_\nu(b)}{2\sin \pi\nu}.\]
This last expression may be simplified by observing that
\begin{eqnarray*}
{}_1F_2\bl(\!\!\begin{array}{c}1\\k+2, k+2+\nu\end{array}\bl|-\chi\br)&=&\sum_{r=0}^\infty\frac{(-\chi)^r}{(k+2+\nu)_r (k+2)_r}\\
&=&(k+1)!\g(k+2+\nu)\sum_{r=0}^\infty \frac{(-\chi)^r}{\g(k\!+\!2\!+\!\nu\!+\!r) (k\!+\!r\!+\!1)!}\\
&=&\frac{(k+1)! \g(k+2+\nu)}{\g(1+\nu)}\,(-\chi)^{-k-1}\sum_{m=k+1}^\infty \frac{(-\chi)^m}{(1+\nu)_m m!}\\
&=&\frac{(k+1)! \g(k+2+\nu)} {(-\chi)^{k+1}}\bl\{\chi^{-\nu/2} J_\nu(b)-\frac{1}{\g(1+\nu)}\sum_{m=0}^k\frac{(-\chi)^m}{(1+\nu)_m m!}\br\}.
\end{eqnarray*}
Hence, after a little simplification, the residue of ${\hat H}(s)$ at $s=1-2k$, $k=1, 2, \ldots$ is given by
\bee\label{e42}
\frac{(-)^k\chi^{\nu/2}}{2} \g(-\nu)F^{(1)}_k(\nu;\chi),
\ee
where $F^{(1)}_k(\nu;\chi)$ is defined in (\ref{e25}).

The singularity of ${\hat H(s)}\zeta(s)a^s$ at $s=1$ is a double pole and consequently we need to expand the quantity in (\ref{e41}) to include the O$(\epsilon)$ term when $s=1+\epsilon$. We find
\[{\hat H}(1+\epsilon)=-\chi^{1+(\nu-\epsilon)/2}\,\frac{\pi\g(-1\!-\!\nu\!+\!\fs\epsilon)}{4\g(2\!-\!\fs\epsilon) \sin \fs\pi\epsilon}\,{}_1F_2\bl(\!\!\begin{array}{c}1\\2\!-\!\fs\epsilon, 2\!+\!\nu\!-\!\fs\epsilon\end{array}\bl|-\chi\br)-\frac{\pi^2J_\nu(b)}{4\sin\pi\nu\,\sin\fs\pi\epsilon}.\]
By Lemma 1 with $\epsilon_1=0$, $\epsilon_2=\epsilon_3=-\fs\epsilon$ it follows that
\[{}_1F_2\bl(\!\!\begin{array}{c}1\\2\!-\!\fs\epsilon, 2\!+\!\nu\!-\!\fs\epsilon\end{array}\bl|-\chi\br)
=\frac{(1+\nu)F(\chi)}{\chi}\bl\{1+\frac{\epsilon}{2} F^{(1)}_*+O(\epsilon^2)\br\},\]
where
\[F(\chi):=\frac{\chi}{1+\nu}\sum_{r=0}^\infty\frac{(-\chi)^r}{(2+\nu)_r (2)_r}=1-\chi^{-\nu/2}\g(1+\nu) J_\nu(b)\]
and
\[F^{(1)}_*:=\frac{\chi}{(1+\nu)F(\chi)}\sum_{r=1}^\infty \frac{(-\chi)^r}{(2+\nu)_r (2)_r}\,\{\Delta\psi(2+r)+\Delta\psi(2+\nu+r)\}\]
\[\hspace{1.2cm}=-\frac{1}{F(\chi)} \sum_{r=1}^\infty\frac{(-\chi)^r}{(1+\nu)_r r!}\,\bl\{\Delta\psi(1+r)+\Delta\psi(1+\nu+r)-\frac{2+\nu}{1+\nu}\br\}\]
upon use of the identity $\psi(1+x)=\psi(x)+1/x$. Then a little algebra produces
\[{\hat H}(1+\epsilon)=\frac{\chi^{\nu/2} \g(-\nu)}{2\epsilon}\bl\{1+\frac{1}{2}\epsilon F(\chi)\{\psi(2)+\psi(-1-\nu)-\log\,\chi+F^{(1)}_*\}+O(\epsilon^2)\br\}.\]
The residue of $H(s)\zeta(s)a^s$ at $s=1$ is therefore
\bee\label{e43}
\mbox{Res}\,(1)=a^{\gamma-1}\bl\{H_2(1,-\nu)+\frac{1}{2}\chi^{\nu/2}\g(-\nu)\bl(\log\,a+\gamma_0-\frac{1}{2}F(\chi)\{\psi(2)+\psi(-1-\nu)-\log\,\chi+F^{(1)}_*\}\br)\br\},
\ee
where 
\[H_2(1,-\nu)=-\bl(\frac{\pi}{2\sin\pi\nu}\br)^2 J_{-\nu}(b).\]

We then have the following result:
\begin{theorem}$\!\!\!.$ \ Let $\gamma+\nu=-1$ with $\nu$ non-integer. Then
\bee\label{e44}
\sum_{n=1}^\infty \frac{n^{-\nu-1} K_\nu(nb/a)}{n^2+a^2}-\mbox{Res}\,(1)\sim {\hat R}(a;\nu)+R(a;-\nu)
\ee
as $|a|\to\infty$ in $|\arg\,a|<\fs\pi$, where $\mbox{Res}\,(1)$ and $R(a;-\nu)$ are defined by (\ref{e43}) and (\ref{e24c}), respectively. The asymptotic sum ${\hat R}(a;\nu)$ is
\[{\hat R}(a;\nu)=\frac{1}{2}a^{\gamma-1}\chi^{\nu/2} \g(-\nu) \sum_{k=1}^\infty (-1)^k \zeta(1-2k) F_k^{(1)}(\nu;\chi)a^{-2k},\]
where $F_k^{(1)}(\nu;\chi)$ is defined in (\ref{e25}) and $\chi=b^2/4$.
\end{theorem}

\noindent{\bf Remark 4.1}\ \
The case $\mu=2$ with $\gamma+\nu=-1$ can be obtained from the identity
\[S_{\nu,\gamma}^{\mu}(a,b)=\frac{1}{2a^2\mu}\bl\{\frac{b}{a} S_{\nu+1,\gamma+1}^{\mu-1}(a,b)-\bl(\nu+a\frac{\partial}{\partial a}\br) S_{\nu,\gamma}^{\mu-1}(a,b)\br\}\]
obtained by differentiation of the series for $S_{\nu,\gamma}^{\mu-1}(a,b)$ together with the relation $K_\nu'(x)=-K_{\nu+1}(x)+(\nu/x) K_\nu(x)$. Then
\[S^2_{\nu,\gamma}(a,b)=\sum_{n=1}^\infty \frac{n^\gamma K_\nu(nb/a)}{(n^2+a^2)^2}=\frac{1}{2a^2}\bl\{\frac{b}{a} S^1_{\nu+1,\gamma+1}(a,b)-\bl(\nu+a \frac{\partial}{\partial a}\br) S^1_{\nu,\gamma}(a,b)\br\},\]
where the asymptotic expansions of $S^1_{\nu,\gamma}(a,b)$ and $S^1_{\nu+1,\gamma+1}(a,b)$ when $\gamma+\nu=-1$ are given by (\ref{e44}) and (\ref{e24}). The derivative of $S^1_{\nu,\gamma}(a,b)$ can be obtained by differentiation of the expansion in (\ref{e44}) to find
\[\frac{\partial}{\partial a} \mbox{Res}\,(1)=\frac{\gamma-1}{a}\,\mbox{Res}\,(1)+\frac{1}{2}a^{\gamma-2}\chi^{\nu/2} \g(-\nu),\]
\[\frac{\partial}{\partial a}{\hat R}(a;\nu)=\frac{1}{2}a^{\gamma-2}\chi^{\nu/2} \g(-\nu) \sum_{k=1}^\infty
(-)^{k+1}(2k+1-\gamma)\zeta(1-2k) F_k^{(1)}(\nu;\chi) a^{-2k},\]
\[\frac{\partial}{\partial a}R(a;-\nu)=\frac{1}{2}a^{\nu-3}\chi^{-\nu/2} \g(\nu) \sum_{k=0}^\infty(-1)^{k+1}(2k+2-\nu)\zeta(-\gamma+\nu-2k) F_k(-\nu) a^{-2k}.\]
\vspace{0.4cm}

\noindent{\bf 4.2\ Integer $\nu$.}\ \ 
When $\nu$ is an integer the singularity at $s=1$ becomes a treble pole. We consider only the case $\nu=0$, $\mu=1$ when $\gamma+\nu=-1$. The poles at $s=1-2k$, $k=1, 2, \ldots$ are double and, from (\ref{e31a}), produce the residue contribution
\bee\label{e45}
a^{-2}\sum_{k=1}^\infty (-)^k\zeta(1-2k)F^{(1)}_k(0;\chi)a^{-2k}\bl\{\log\,\frac{2a}{b}-\gamma_0+\frac{\zeta'(1-2k)}{\zeta(1-2k)}+{\hat F}^{(1)}_*\br\},
\ee
where $F^{(1)}_k(0;\chi)$ is given by (\ref{e25}) and\footnote{In this example there are two instances of the quantity $F^{(1)}_*$. To distinguish between them we add a circumflex over the first occurring quantity.}
\bee\label{e45a}
{\hat F}^{(1)}_*=\frac{1}{F^{(1)}_k(0;\chi)}\sum_{r=1}^k\frac{(-\chi)^r}{(r!)^2}\,\Delta\psi(1+r).
\ee

The functional form of $H(s)$ in this case is \cite{M}
\[H(s)=\frac{1}{4}\chi^{3/2-s/2}\g^2(\fs s-\f{3}{2})\,{}_1F_2\bl(\!\!\begin{array}{c}1\\ \f{5}{2}\!-\!\fs s, \f{5}{2}\!-\!\fs s\end{array}\bl|\,-\chi\br)\hspace{6cm}\]
\[\hspace{5cm}+\frac{\pi^2}{4\sin^2 \fs\pi(s\!-\!1)}\bl\{J_0(b) \cos \fs\pi(s\!-\!1)-Y_0(b) \sin \fs\pi(s\!-\!1)\br\},\]
where $Y_0(b)$ is the Bessel function of the second kind of zero order. Then, with $s=1+\epsilon$, $\epsilon\to0$, 
\[H(1+\epsilon)=\frac{\Lambda(\epsilon)}{\epsilon^2}\,{}_1F_2\bl(\!\!\begin{array}{c}1\\ 2\!-\!\fs \epsilon, 2\!-\!\fs \epsilon\end{array}\bl|\,-\chi\br)+\frac{1}{\epsilon^2}\bl\{J_0(b)-\frac{1}{2}\pi\epsilon Y_0(b)-\frac{\pi^2\epsilon^2}{24} J_0(b)+O(\epsilon^3)\br\},\]
where
\[\Lambda(\epsilon)=\frac{\chi^{-\epsilon/2} \g^2(1+\fs\epsilon)}{(1-\fs\epsilon)^2}=\sum_{r=0}^\infty A_r \epsilon^r\]
with 
\[A_0=1,\quad A_1=1-\gamma_0-\frac{1}{2}\log\,\chi,\] \[A_2=\frac{3}{4}-\gamma_0+\frac{1}{2}\gamma_0^2+\frac{\pi^2}{24}+\frac{1}{2}\log\,\chi\bl(\gamma_0-1+\frac{1}{4}\log\,\chi\br).\]

From (\ref{e16}) we obtain the expansion
\[{}_1F_2\bl(\!\!\begin{array}{c}1\\ 2\!-\!\fs \epsilon, 2\!-\!\fs \epsilon\end{array}\bl|\,-\chi\br)=
F(\chi)\{1+\epsilon F^{(1)}_*+\epsilon^2 F^{(2)}_*+O(\epsilon^3)\},\]
where
\[F(\chi)= \sum_{r=0}^\infty \frac{(-\chi)^r}{((r+1)!)^2}=\frac{1-J_0(b)}{\chi}\]
and
\[F^{(1)}_*=\frac{1}{F(\chi)}\sum_{r=1}^\infty \frac{(-\chi)^r}{((r+1)!)^2}\,\Delta\psi(2+r),\quad F^{(2)}_*=\frac{1}{F(\chi)}\sum_{r=1}^\infty\frac{(-\chi)^r}{((r+1)!)^2}\,\{2\Delta^2\psi(2+r)-\Delta\psi'(2+r)\}.\]
The expansion of $\zeta(s)a^s$ about $s=1$ follows from (\ref{e15}) in the form
\[\zeta(1+\epsilon)a^{1+\epsilon}=\frac{a}{\epsilon}\sum_{r=0}^\infty B_r\epsilon^r,\]
where
\[B_0=1,\quad B_1=\log\,a+\gamma_0, \quad B_2=\frac{1}{2}(\log\,a)^2+\gamma_0 \log\,a-\gamma_1.\]
Then, after some straightforward algebra, we find
\[\mbox{Res}\,(1)=a^{-2}\bl\{B_2+(1-J_0(b))\{(A_1+A_0F^{(1)}_*)B_1+(A_2+A_1 F^{(1)}_*+A_0 F^{(2)}_*)B_0\}\]
\bee\label{e46}
\hspace{5cm}-\frac{\pi^2}{24}B_0J_0(b)-\frac{1}{2}\pi B_1Y_0(b)\br\}.
\ee

Combination of (\ref{e45}) and (\ref{e46}) then yields the following result:
\begin{theorem}$\!\!\!.$ \ When $\nu=0$, $\gamma=-1$ and $\mu=1$ we have the expansion
\bee\label{e47}
\sum_{n=1}^\infty\frac{K_0(nb/a)}{n(n^2+a^2)}-\mbox{Res}\,(1)\sim
a^{-2}\sum_{k=1}^\infty (-)^k\zeta(1-2k)F^{(1)}_k(0;\chi)a^{-2k}\bl\{\log\,\frac{2a}{b}-\gamma_0+\frac{\zeta'(1-2k)}{\zeta(1-2k)}+{\hat F}^{(1)}_*\br\}
\ee
as $|a|\to\infty$ in $|\arg\,a|<\fs\pi$, where ${\hat F}^{(1)}_*$ and $\mbox{Res}\,(1)$  are defined in (\ref{e45a}) and (\ref{e46}). 
\end{theorem}

\vspace{0.6cm}

\begin{center}
{\bf 5. \ Numerical results and concluding remarks}
\end{center}
\setcounter{section}{5}
\setcounter{equation}{0}
\renewcommand{\theequation}{\arabic{section}.\arabic{equation}}
We present some numerical results to illustrate the accuracy of the various expansions obtained. The values shown give the absolute relative error in the high-precision computation of the Mathieu-Bessel series from (\ref{e11}) using the asymptotic expansions in Theorems 1--6. In each case the asymptotic sums are optimally truncated (that is, truncation just before the least term in the expansion). 

In Table 1    
we show values\footnote{In the tables we write the values as $x(y)$ instead of $x\times 10^y$.} of $S_{\nu,\gamma}^\mu(a,b)$
and the relative error for different values of $a$ and $b$ when Theorem 1 applies. Table 2 gives examples of the relative error for complex values of $a$. It is seen that the error progressively decreases as $\arg\,a$ approaches $\fs\pi$.
Table 3 shows examples for different $a$ with $\nu$, $\gamma$ and $\mu$ chosen to correspond to Theorems 2--6. 

\begin{table}[t]
\caption{\footnotesize{The absolute relative error in the computation of $S^\mu_{\nu,\gamma}(a,b)$ from Theorem 1 for different $a$ and $b$ when $\mu=1$, $\nu=2/3$ and $\gamma=1$. The optimal truncation index $k_o$ employed in the asymptotic sums $R(a,\pm\nu)$ is indicated.}}
\begin{center}
\begin{tabular}{|rr|cc|cc|}
\hline
&&&&&\\[-0.25cm]
\mcol{2}{|c|}{} & \mcol{2}{c|}{$b=1$} &\mcol{2}{c|}{$b=2$} \\
\mcol{1}{|c}{$a$} & \mcol{1}{c|}{$k_o$} & \mcol{1}{c}{$S^\mu_{\nu,\gamma}(a,b)$} & \mcol{1}{c|}{Error}& \mcol{1}{c}{$S^\mu_{\nu,\gamma}(a,b)$} & \mcol{1}{c|}{Error}\\
[.1cm]\hline
&&&&&\\[-0.25cm]
2 & 5 &  $4.67988(-01)$ & $1.050(-04)$ & $1.41942(-01)$ & $3.177(-05)$ \\
4 & 11 & $5.36953(-01)$ & $9.670(-10)$ & $1.82047(-01)$ & $3.964(-10)$ \\
6 & 17 & $5.56269(-01)$ & $5.951(-15)$ & $1.93851(-01)$ & $2.807(-15)$ \\
8 & 24 & $5.64899(-01)$ & $2.899(-20)$ & $1.99192(-01)$ & $8.323(-21)$ \\
10 & 30 & $5.69662(-01)$ & $1.394(-25)$ & $2.02157(-01)$ & $4.939(-26)$\\
[.1cm]\hline
\end{tabular}
\end{center}
\end{table}
\begin{table}[h]
\caption{\footnotesize{The absolute relative error in the computation of $S^\mu_{\nu,\gamma}(a,b)$ using the expansion in Theorem 1 at optimal truncation when $a=8\exp (i\theta)$ as a function $\theta$  when $\gamma=1$, $\nu=2/3$ and $b=1$.}}
\begin{center}
\begin{tabular}{|l|c|c|c|}
\hline
&&&\\[-0.25cm]
\mcol{1}{|c|}{$\theta$} &\mcol{1}{c|}{$\mu=1/2$} &  \mcol{1}{c|}{$\mu=1$} &\mcol{1}{c|}{$\mu=5/2$} \\
[.1cm]\hline
&&&\\[-0.25cm]
0        & $3.323(-21)$ &  $2.899(-20)$ & $2.689(-18)$  \\
$\pi/8$  & $1.673(-19)$ &  $1.478(-18)$ & $1.417(-16)$  \\
$\pi/4$  & $9.057(-15)$ &  $7.982(-14)$ & $7.544(-12)$  \\
$3\pi/8$ & $1.094(-07)$ &  $9.627(-07)$ & $8.977(-05)$  \\
[.1cm]\hline
\end{tabular}
\end{center}
\end{table}
\begin{table}[h]
\caption{\footnotesize{The absolute relative error in the computation of $S^\mu_{\nu,\gamma}(a,b)$ using the expansions in Theorems 2--6 at optimal truncation for different $a$ when $b=1$.}}
\begin{center}
\begin{tabular}{|r|c|c|c|c|c|}
\hline
&&&&&\\[-0.25cm]
\mcol{1}{|c|}{} &\mcol{1}{|c|}{$\mu=2$} &\mcol{1}{c|}{$\mu=2$} &  \mcol{1}{c|}{$\mu=\f{8}{5}$} &\mcol{1}{c|}{$\mu=1$} &\mcol{1}{c|}{$\mu=1$}\\
\mcol{1}{|c|}{$a$} &\mcol{1}{|c|}{$\nu=0,\ \gamma=1$} &\mcol{1}{c|}{$\nu=1,\ \gamma=\f{3}{4}$} &  \mcol{1}{c|}{$\nu=\f{1}{3},\ \gamma=-\f{4}{3}$} &\mcol{1}{c|}{$\nu=\f{1}{4}, \gamma=-\f{5}{4}$} &\mcol{1}{c|}{$\nu=0, \gamma=-1$}\\
[.1cm]\hline
&&&&&\\[-0.25cm]
2 & $1.204(-03)$ & $1.039(-04)$ &  $7.667(-06)$ & $9.376(-06)$     & $1.324(-05)$\\
4 & $2.857(-08)$ & $1.398(-09)$ &  $2.539(-11)$ & $1.581(-11)$ & $2.080(-11)$\\
6 & $3.095(-13)$ & $1.081(-14)$ &  $8.702(-17)$ & $3.545(-17)$ & $5.127(-17)$\\
8 & $2.366(-18)$ & $6.641(-20)$ &  $2.563(-22)$ & $9.315(-23)$ & $1.417(-22)$\\
10& $1.554(-23)$ & $2.802(-25)$ &  $9.106(-28)$ & $2.420(-28)$ & $4.269(-28)$\\
[.1cm]\hline
\end{tabular}
\end{center}
\end{table}

We remark that the asymptotic expansion of the alternating version of (\ref{e11}) can be deduced by making use of the identity
\bee\label{e51}
\sum_{n=1}^\infty \frac{(-)^{n-1} n^\gamma }{(n^2+a^2)^\mu}\,K_\nu(nb/a)=S_{\nu,\gamma}^\mu(a,b)-2^{1+\gamma-2\mu} S_{\nu,\gamma}^\mu(\fs a,b).
\ee
From Theorem 1, this then produces the asymptotic expansion (when $\gamma+\nu\neq -1, -3, \ldots$)
\bee\label{e52}\sum_{n=1}^\infty \frac{(-)^{n-1} n^\gamma }{(n^2+a^2)^\mu}\,K_\nu(nb/a)\sim
{\tilde R}(a;\nu)+{\tilde R}(a;-\nu)
\ee
as $|a|\to\infty$ in $|\arg\,a|<\fs\pi$, where
\[{\tilde R}(a;\nu):=\frac{1}{2}a^{-\nu-2\mu}\chi^{\nu/2}\g(-\nu) \sum_{k=0}^\infty \frac{(-)^k(\mu)_k}{k!}\{1-2^{1+\omega_k}\}\zeta(-\omega_k)F_k^{(\mu)}(\nu;\chi) a^{-2k}\]
and $\omega_k=\gamma+\nu+2k$.

Finally, we observe that the asymptotic expansion for large $|a|$ of the series involving the derivative of the Bessel function can be obtained from
\[\sum_{n=1}^\infty\frac{n^\gamma K'_\nu(nb/a)}{(n^2+a^2)^\mu}=-\frac{1}{2}\bl\{S^\mu_{\nu-1,\gamma}(a,b)+S^\mu_{\nu+1,\gamma}(a,b)\br\},\]
which follows from the relation $K_\nu(x)=-\fs\{K_{\nu-1}(x)+K_{\nu+1}(x)\}$.

\vspace{0.6cm}

\begin{center}
{\bf Appendix A: Discussion of the pole structure of $H(s)$}
\end{center}
\setcounter{section}{1}
\setcounter{equation}{0}
\renewcommand{\theequation}{\Alph{section}.\arabic{equation}}\vspace{0.6cm}
We examine the pole structure of the function $H(s)$ defined in (\ref{e22a}), which has apparent singularities when $\lambda_{\pm\nu}-\mu=k$, $k=0, \pm1, \pm2, \ldots\,$, where $\lambda_{\nu}=\fs(\gamma+\nu+s)$; that is, at $s=2\mu-\gamma\mp\nu+2k$. It will be shown that $H(s)$ is in fact regular at these points due to a cancellation of terms.

It will be sufficient to deal with the situation corresponding to $\lambda_\nu-\mu=k$, since the case with $\nu\to-\nu$ is similar. The functions $H_1(s,\nu)$ and $H_2(s,\nu)$ appearing in (\ref{e22a}) present singularities at these points. 
To demonstrate that the combination $H_1(s,\nu)+H_2(s,\nu)$ (and hence $H(s)$) is regular there, it will be sufficient to show that the function
\[Q(s):=\frac{4}{\pi^2} \sin \pi(\lambda_{\nu}\!-\!\mu)\{H_1(s,\nu)+H_2(s,\nu)\}\]
has a simple zero at these points. 

We first consider the points $s_k=2\mu-\gamma-\nu+2k$,  with $k=0, 1, 2, \ldots\,$ to obtain
\bee\label{a1}
Q(s_k)=\frac{\chi^{\nu/2-k}}{\sin \pi(k\!-\!\nu)}{}_1{\bf F}_2\bl(\!\!\begin{array}{c}\mu\\1-k, 1+\nu-k\end{array}\bl|\,-\chi\br)+\frac{\chi^{\nu/2} \g(\mu+k)}{\g(\mu) \sin \pi\nu}\,
{}_1{\bf F}_2\bl(\!\!\begin{array}{c}\mu+k\\1+k, 1+\nu\end{array}\bl|\,-\chi\br).
\ee

The first term on the right-hand side of (\ref{a1}) can be written as
\[\frac{\chi^{\nu/2-k}}{\sin \pi(k\!-\!\nu)}{}_1{\bf F}_2\bl(\!\!\begin{array}{c}\mu\\1-k, 1+\nu-k\end{array}\bl|\,-\chi\br)
=\frac{(-)^{k+1}\chi^{\nu/2-k}}{\sin \pi\nu\,\g(1\!+\!\nu\!-\!k)}\sum_{r=k}^\infty\frac{(\mu)_r (-\chi)^r}{\g(1\!-\!k\!+\!r)(1\!+\!\nu\!-\!k)_r r!}\]
\[=-\frac{\chi^{\nu/2}}{\sin \pi\nu\,\g(1\!+\!\nu\!-\!k)}\sum_{n=0}^\infty\frac{(\mu)_{n+k} (-\chi)^n}{(n+k)! (1+\nu-k)_{n+k} n!}=-\frac{\chi^{\nu/2} \g(\mu+k)}{\g(\mu) \sin \pi\nu\,k! \g(1+\nu)} \sum_{n=0}^\infty\frac{(\mu+k)_n (-\chi)^n}{(1+k)_n(1+\nu)_n  n!}\]
\[=-\frac{\chi^{\nu/2} \g(\mu+k)}{\g(\mu) \sin \pi\nu}\,
{}_1{\bf F}_2\bl(\!\!\begin{array}{c}\mu+k\\1+k, 1+\nu\end{array}\bl|\,-\chi\br).\]
Hence $Q(s_k)=0$ when $k=0, 1, 2, \ldots\,$.

Next, consider the points $s_{k}$ with $k=-1, -2, \ldots\ $. Then we have
\bee\label{a2}
Q(s_{k})=\frac{(-)^{k+1}\chi^{\nu/2+k}}{\sin \pi\nu}{}_1{\bf F}_2\bl(\!\!\begin{array}{c}\mu\\1+k,1+\nu+k\end{array}\bl|\,-\chi\br)+\frac{\chi^{\nu/2} \g(\mu-k)}{\g(\mu) \sin \pi\nu}{}_1{\bf F}_2\bl(\!\!\begin{array}{c}\mu-k\\1-k, 1+\nu\end{array}\bl|\,-\chi\br).
\ee
The second term in (\ref{a2}) can be written as
\[\frac{\chi^{\nu/2} \g(\mu-k)}{\g(\mu) \sin \pi\nu}{}_1{\bf F}_2\bl(\!\!\begin{array}{c}\mu-k\\1-k, 1+\nu\end{array}\bl|\,-\chi\br)=\frac{\chi^{\nu/2} \g(\mu-k)}{\g(\mu) \g(1+\nu) \sin \pi\nu}
\sum_{r=k}^\infty\frac{(\mu-k)_r (-\chi)^r}{\g(1-k+r) (1+\nu)_r r!}\]
\[=\frac{(-)^k \chi^{\nu/2+k}\g(\mu-k)}{\g(\mu)\g(1+\nu)\sin \pi\nu}\sum_{n=0}^\infty\frac{(\mu-k)_{n+k} (-\chi)^n}{(n+k)! (1+\nu)_{n+k} n!}=\frac{(-)^k \chi^{\nu/2+k}}{\sin \pi\nu\,\g(1+\nu+k) k!}\sum_{n=0}^\infty\frac{(\mu)_n (-\chi)^n}{(1+k)_n(1+\nu+k)_n n!}\]
\[=\frac{(-)^{k}\chi^{\nu/2+k}}{\sin \pi\nu}{}_1{\bf F}_2\bl(\!\!\begin{array}{c}\mu\\1+k,1+\nu+k\end{array}\bl|\,-\chi\br).\]
Hence $Q(s_{k})=0$ when $k=-1, -2, \dots\,$.

Thus $H(s)$ is regular at the points $s_k=2\mu-\gamma-\nu+2k$, $k=0, \pm 1 \pm2,\ldots\,$. A similar argument applies at the points $s_k=2\mu-\gamma+\nu+2k$ (by putting $\nu\to -\nu$). This concludes the demonstration that $H(s)$ is regular at the points $\lambda_{\pm\nu}-\mu$, $k=0, \pm1, \pm2, \ldots\,$.

\vspace{0.6cm}

\begin{center}
{\bf Appendix B: Evaluation of the residues when $\nu=1$}
\end{center}
\setcounter{section}{2}
\setcounter{equation}{0}
\renewcommand{\theequation}{\Alph{section}.\arabic{equation}}
We require to determine the residue contribution in Theorem 1 given by
\[R(a;\nu)+R(a;-\nu)\]
when $\nu=1$. Let $\nu=1+\epsilon$, $\epsilon\to0$ so that, from (\ref{e24c}), we have
\[R(a;1+\epsilon)=\frac{a^{-1-2\mu-\epsilon} \chi^{(1+\epsilon)/2}}{2\epsilon \g(2+\epsilon)} \sum_{k=0}^\infty \frac{(-)^k}{k!} (\mu)_k \zeta(-\gamma\!-\!1\!-\!2k\!-\!\epsilon)a^{-2k} F_k^{(\mu)}(1+\epsilon;\chi) +O(\epsilon^2)\]
and
\[R(a;-1-\epsilon)=\frac{1}{2}a^{1-2\mu+\epsilon} \chi^{-(1+\epsilon)/2}\g(1+\epsilon)\sum_{k=0}^\infty \frac{(-)^k}{k!} (\mu)_k\zeta(1\!-\!\gamma\!-\!2k\!+\!\epsilon)a^{-2k} F_k^{(\mu)}(-1-\epsilon;\chi).\]

A straightforward rearrangement of the series expansion in (\ref{e24a}) shows that
\[F_k^{(\mu)}(-1-\epsilon:\chi)=\sum_{r=0}^k\frac{(-k)_r (-\chi)^r}{(-\epsilon)_r(1-\mu-k)_r r!}=1+\frac{k\chi\sigma_{k-1}}{\epsilon(k+\mu-1)},\]
where
\[\sigma_k:=\sum_{r=0}^{k}\frac{(-k)_r (-\chi)^r}{(1-\epsilon)_r(1-\mu-k)_r(2)_r}.\]
Thus we obtain
\[R(a;-1-\epsilon)=\frac{1}{2}a^{1-2\mu} \chi^{-1/2}\sum_{k=0}^\infty \frac{(-)^k}{k!} (\mu)_k\zeta(1\!-\!\gamma\!-\!2k)a^{-2k}+O(\epsilon)\hspace{5cm}\]
\bee\label{b1}
\hspace{4cm}+a^{1-2\mu+\epsilon} \chi^{(1-\epsilon)/2}\frac{\g(1+\epsilon)}{2\epsilon}\sum_{k=1}^\infty \frac{(-)^k}{(k-1)!} (\mu)_{k-1}\zeta(1\!-\!\gamma\!-\!2k\!+\!\epsilon)a^{-2k}\,\sigma_{k-1}.
\ee
The second term on the right-hand side of (\ref{b1}) becomes, upon putting $k\to k+1$,
\[
-a^{-1-2\mu+\epsilon} \chi^{(1-\epsilon)/2}\frac{\g(1+\epsilon)}{2\epsilon}\sum_{k=0}^\infty \frac{(-)^k}{k!} (\mu)_{k}\zeta(-1\!-\!\gamma\!-\!2k\!+\!\epsilon)a^{-2k}\,\sigma_{k}\]
\[=-\frac{a^{-1-2\mu}}{2\epsilon} \chi^{1/2} \sum_{k=0}^\infty \frac{(-)^k}{k!} (\mu)_k \zeta(-1\!-\!\gamma\!-\!2k) a^{-2k} F_k^{(\mu)}(1;\chi)\bl\{1+A\epsilon+O(\epsilon^2)\br\},\]
where $F_k^{(\mu)}(1;\chi)$ is defined in (\ref{e24a}) and, using an obvious extension of Lemma 1,
\[A=\log\,\frac{2a}{b}-\gamma_0+\frac{\zeta'(-1\!-\!\gamma\!-\!2k)}{\zeta(-1\!-\!\gamma\!-\!2k)}+\sum_{r=0}^k\frac{(-k)_r (-\chi)^r}{(2)_r (1-\mu-k)_r r!}\,\Delta\psi(1+r).\]\

Application of Lemma 1 shows that
\bee\label{b2}
R(a;1+\epsilon)=\frac{a^{-1-2\mu}}{2\epsilon} \chi^{1/2} \sum_{k=0}^\infty\frac{(-)^k}{k!} (\mu)_k \zeta(-1\!-\!\gamma\!-\!2k) a^{-2k} F_k^{(\mu)}(1;\chi)\bl\{1+B\epsilon+O(\epsilon^2)\br\},
\ee
where
\[B=-\log\,\frac{2a}{b}-\psi(2)-\frac{\zeta'(-1\!-\!\gamma\!-\!2k)}{\zeta(-1\!-\!\gamma\!-\!2k)}-\sum_{r=0}^k\frac{(-k)_r (-\chi)^r}{(2)_r (1-\mu-k)_r r!}\,\Delta\psi(2+r).\]
Combination of (\ref{b1}) and (\ref{b2}) then produces the limiting value
\[R(a;1+\epsilon)+R(a;-1-\epsilon)=\frac{1}{2}a^{1-2\mu} \chi^{-1/2}\sum_{k=0}^\infty \frac{(-)^k}{k!} (\mu)_k\zeta(1\!-\gamma\!-\!2k)a^{-2k}\]
%\[+a^{-1-2\mu}\chi^{1/2} \sum_{k=0}^\infty\frac{(-)^{k+1}}{k!} (\mu)_k \zeta(-1\!-\!\gamma\!-\!2k) a^{-2k} F_k^{(\mu)}(1;\chi)\]
\bee\label{b3}
+a^{-1-2\mu}\chi^{1/2}\sum_{k=0}^\infty\frac{(-)^{k+1}}{k!} (\mu)_k \zeta(-\gamma-1-2k) F_k^{(\mu)}(1;\chi) a^{-2k}
\bl\{\log\,\frac{2a}{b}-\gamma_0+\frac{1}{2}+\frac{\zeta'(-\gamma\!-\!1\!-\!2k)}{\zeta(-\gamma\!-\!1\!-\!2k)}+F_*^{(1)}\br\}
\ee
as $\epsilon\to0$, where
\[F_*^{(1)}:=\frac{1}{2F_k^{(\mu)}(1;\chi)}\sum_{r=0}^k\frac{(-k)_r (-\chi)^r}{(2)_r(1-\mu-k)_r r!} \{\Delta\psi(1+r)+\Delta(2+r)\}\]
\bee\label{b4}
\hspace{0.9cm}=\frac{1}{F_k^{(\mu)}(1;\chi)}\sum_{r=0}^k\frac{(-k)_r (-\chi)^r}{(2)_r(1-\mu-k)_r r!} \bl\{\Delta\psi(1+r)-\frac{r}{2(1+r)}\br\}.
\ee
\vspace{0.6cm}

\begin{center}
{\bf Appendix C: Expansion of ${}_1F_2(-\chi)$ appearing in (\ref{e36a}) }
\end{center}
\setcounter{section}{3}
\setcounter{equation}{0}
\renewcommand{\theequation}{\Alph{section}.\arabic{equation}}
We require the expansion of
\[{\cal F}\equiv{}_1F_2\bl(\!\!\begin{array}{c}-n+1+\epsilon\\ 1-\alpha-n+\epsilon, 1+\nu\end{array}\bl|\,-\chi\br)
=\bl\{\sum_{r=0}^{n-1}+\sum_{r=n}^\infty\br\}\frac{(-n+1+\epsilon)_r (-\chi)^r}{(1+\nu)_r(1-\alpha-n+\epsilon)_r r!}\]
as $\epsilon\to0$ when $n$ is a positive integer and $\mu\neq 1, 2, \ldots\,$. In (\ref{e36a}) we have $n=m+1$ and $\alpha=\mu-1$.

Now, from (\ref{e16}) and use of the standard properties of the $\psi$ function
\[(-n+1+\epsilon)_r=(-n+1)_r\{1+\epsilon \Delta\psi(n\!-\!r)+O(\epsilon^2)\}\qquad (0\leq r\leq n-1),\]
so that
\[\sum_{r=0}^{n-1}\frac{(-n+1+\epsilon)_r (-\chi)^r}{(1+\nu)_r(1-\alpha-n+\epsilon)_r r!}\hspace{6cm}\]
\bee\label{c1}
=F_{n-1}^{(\mu)}(\nu;\chi)+\epsilon
\sum_{r=0}^{n-1}\frac{(-n+1)_r (-\chi)^r}{(1+\nu)_r(1-\alpha-n)_r r!}\,(\Delta\psi(n\!-\!r)-\Delta\psi(n\!+\!\alpha\!-\!r))+O(\epsilon^2).
\ee
The infinite sum can be rearranged in the form
\[\sum_{r=n}^\infty\frac{(-n+1+\epsilon)_r (-\chi)^r}{(1+\nu)_r(1\!-\!\alpha\!-\!n\!+\!\epsilon)_r r!}
=\frac{(-n+1+\epsilon)_n (-\chi)^n}{(1+\nu)_n(1\!-\!\alpha\!-\!n\!+\!\epsilon)_n n!}\,{}_2F_3\bl(\!\!\begin{array}{c}1, 1+\epsilon\\n+1, n+1+\nu, 1-\alpha+\epsilon\end{array}\bl|\,-\chi\br)\]
\bee\label{c2}
=-\frac{\epsilon (-\chi)^{n}}{n(1+\nu)_{n}(\alpha)_{n}}\,{}_2F_3\bl(\!\!\begin{array}{c} 1,1\\n+1, n+1+\nu, 1-\alpha\end{array}\bl|\,-\chi\br)+O(\epsilon^2).
\ee

Then combination of (\ref{c1}) and (\ref{c2}) yields
\bee\label{c3}
{\cal F}=F_{n-1}^{(\mu)}(\nu;\chi)\{1+\epsilon F_*^{(1)}+O(\epsilon^2)\},
\ee 
where
\[F_*^{(1)}=\frac{1}{F_{n-1}^{(\mu)}(\nu;\chi)}\bl\{\sum_{r=0}^{n-1}\frac{(-n+1)_r (-\chi)^r}{(1+\nu)_r(1-\alpha-n)_r r!}\,(\Delta\psi(n\!-\!r)-\Delta\psi(n\!+\!\alpha\!-\!r))\hspace{4cm}\]
\bee\label{c4}
\hspace{4cm}-\frac{\epsilon (-\chi)^{n}}{n(1+\nu)_{n}(\alpha)_{n}}\,{}_2F_3\bl(\!\!\begin{array}{c} 1,1\\n+1, n+1+\nu, 1-\alpha\end{array}\bl|\,-\chi\br)\br\}.
\ee

\end{document}